\documentclass{article}
\usepackage{url}
\usepackage{cite}
\usepackage{amsmath,amssymb,amsfonts}
\usepackage{algorithmic}
\usepackage{graphicx}
\usepackage{textcomp}
\usepackage{xcolor}
\begin{document}

\title{About the "accurate mode" of the IEEE 1788-2015 standard for interval arithmetic
}

\author{Nathalie Revol}
\date{\small
INRIA Lyon - AriC team\\
LIP (UMR 5668, ENS Lyon, University Lyon 1, Inria, CNRS) \\
Lyon, France \\
ORCID 0000-0002-2503-2274}

\maketitle

\begin{abstract}
The IEEE 1788-2015 standard for interval arithmetic defines three accuracy modes for the so-called set-based flavor: tightest, accurate and valid.
This work in progress focuses on the accurate mode.
First, an introduction to interval arithmetic and to the IEEE 1788-2015 standard is given, then the accurate mode is defined.
How can this accurate mode be tested, when a library implementing interval arithmetic claims to provide this mode? The chosen approach is unit testing, and the elaboration of testing pairs for this approach is developed.
A discussion closes this paper: how can the tester be tested? And if we go to the roots of the subject, is the accurate mode really relevant or should it be dropped off in the next version of the standard?
\end{abstract}

\noindent {\bf Keywords.}
interval arithmetic, IEEE 1788-2015 standard, accuracy, accurate mode

\section{Introduction}
The core of this article is a discussion of the "accurate mode" of the IEEE 1788-2015 standard.
Before digging deep in the details of this mode, let us first briefly introduce this standard in general in Section~\ref{subsec:1788}, and the accurate mode in particular in Section~\ref{subsec:accurate-mode-def}.
In Section~\ref{sec:testing}, we will detail our method for testing compliance of interval arithmetic libraries regarding the accurate mode.
This testing experiment led us to question whether keeping the accurate mode in future revisions of the standard is really relevant, as developed in Section~\ref{sec:ccl}.

\subsection{Introduction to interval arithmetic and to the IEEE 1788-2015 standard}
\label{subsec:1788}
The main motivation for interval arithmetic is to handle inaccuracies and uncertainties in data and during calculations.
It operates on intervals rather than on numbers, and the result of each operation is an interval that is guaranteed to contain the range of this operation on the input intervals.
Indeed, the main strength of interval arithmetic is to offer guarantees on the results of computations:
the result contains the exact result if there is one, or the range of all possible results given the uncertainties on the input data.

Let us denote intervals in boldface, for instance \(\mathbf{x}\), and the set of real intervals \( [ \underline{x}, \bar{x}] = \{ x \in \mathbb{R} \: : \: \underline{x} \leq x \leq \bar{x}\}\) for \( \underline{x} \in \mathbb{R} \cup \{ - \infty\}\) and \( \bar{x} \in \mathbb{R} \cup \{+\infty\}\), as \( \mathrm{I}\!\mathbb{R}\).
For any \( \mathbf{x} \in \mathrm{I}\!\mathbb{R}\), for any \( \mathbf{y} \in \mathrm{I}\!\mathbb{R}\), the result of the operation \(\diamond\) on \(\mathbf{x}\) and \(\mathbf{y}\) is the tightest (for inclusion) interval
\(\mathbf{z}\) that contains \(\{ x \diamond y, \: x \in \mathbf{x} \mbox{ and } y \in \mathbf{y} \} \) if
\( (\mathbf{x}, \mathbf{y}) \subset \mbox{Dom}_{\diamond}\).

For most operations, this translates into formulas that are easy to implement, such as
\([\underline{x}, \bar{x}] - [\underline{y},\bar{y}] = [\underline{x}-\bar{y}, \bar{x}-\underline{y}]\).
For non-monotonic operations, more care is required: in what follows, we assume that developers implement such functions with greatest care.
However, in the first variants of interval arithmetic, there was no consensus on cornercases such as 
\(\sqrt{[-1,2]}\).
The interval arithmetic community has gathered efforts and the result is the IEEE 1788-2015 standard for interval arithmetic~\cite{IEEE-1788-2015}.

\subsection{Accuracy of the computed results}
In this work, we will focus only on the accuracy (in the usual sense of this word) of the results, for an operation evaluation by a library under consideration.
More specifically, by accuracy of a result, we mean the difference between the exact, tightest possible result \(\mathbf{r}^*\) and the result \(\mathbf{r}\) computed by the library.
This result \(\mathbf{r}\) must satisfy the {\em fundamental principle of interval arithmetic}: \( \mathbf{r} \supset \mathbf{r}^*\).

However, \(\mathbf{r}\) can be strictly larger than \(\mathbf{r}^*\), in particular if floating-point arithmetic is used as the underlying arithmetic.
For instance, the formula above for \([\underline{x}, \bar{x}] - [\underline{y},\bar{y}] \)
becomes \( \mathbf{r} = [\mathrm{RD}(\underline{x}-\bar{y}), \mathrm{RU}(\bar{x}-\underline{y})]\),
where \(\mathrm{RD}\) stands for ``rounding downwards'' and \(\mathrm{RU}\) for ``rounding upwards'',
whereas \( \mathbf{r}^* = [\underline{x}-\bar{y}, \bar{x}-\underline{y}]\) in exact arithmetic.
In this example,  \( \mathbf{r}\) is the tightest possible result, given the underlying arithmetic.
But \(\mathbf{r}' = [ \underline{r}', \bar{r}']\) for any \( \underline{r}' \leq \mathrm{RD} (\underline{x}-\bar{y}) \) and any \( \bar{r}' \geq \mathrm{RU}(\bar{x}-\underline{y}) \) is also an admissible result.
Even \( \mathbb{R} = ( -\infty, +\infty)\) satisfies the enclosure requirement, however it does not convey useful information.

\subsection{The three accuracy modes of the IEEE 1788-2015 standard}
The IEEE 1788-2015 standard for interval arithmetic handles this issue by defining several possible {\em accuracy modes}, when the result exists. 
Two general accuracy modes are the {\em tightest} mode and the {\em valid} mode. Here are their definitions, given in Section 7.5.4 of~\cite{IEEE-1788-2015}.
``{\em Two accuracy modes for an interval-valued operation are defined for all flavors.}
(Note: A flavor is a variant of interval arithmetic, based on a specific mathematical theory: such theories can be Kaucher arithmetic, modal arithmetic, set theory\ldots). 
\([\ldots]\)
{\em The basic property of enclosure in the sense of the flavor,} \([\ldots]\) 
{\em is called \underline{valid} accuracy mode. Unsurprisingly, the property that the result equals the} tightest representable (given the underlying arithmetic format) hull of the exact mathematical {\em value is called \underline{tightest} accuracy mode.}''

For the set-based flavor, which is the only flavor defined in the current version of the standard, an additional accuracy mode is added, thus three accuracy modes are defined: tightest, accurate and valid. These modes are  ``{\em linearly ordered by strength, with tightest the strongest and valid the weakest, where mode \(M\) is stronger than mode \(M'\) if \(M\) implies \(M'\)} \cite[Section 7.5.4]{IEEE-1788-2015}.

\noindent {\bf Caveat.} Let us insist on a subtlety on the wording: three {\em accuracy} modes exist in the set-based flavor of the IEEE 1788-2015 standard, that characterize how close or far the computed result is from the mathematically exact result. We use the noun {\em accuracy} to indicate these modes. These three modes are: {\em tightest}, {\em accurate} and {\em valid}: here {\em tightest} indicates that the result is the tightest possible one, {\em valid} indicates that the inclusion property is satisfied, and {\em accurate}, which lies in between, is one of the three {\em accuracy} mode. When we use the noun {\em accuracy}, we refer to this set of modes. When we use the adjective {\em accurate}, we refer to one specific mode.

To claim conformance with the IEEE 1788-2015 standard, any implementation of interval arithmetic that implements the set-based flavor must specify 
``{\em What accuracy is achieved (i.e., tightest, accurate, or valid) for each of the implementation's interval operations?}''
(cf. \cite[Section 12.10.3]{IEEE-1788-2015}.

\subsection{The ``accurate mode'' of the IEEE 1788-2015 standard}
\label{subsec:accurate-mode-def}
Here is the definition of the accurate mode, as given in~\cite[Section 12.10.1]{IEEE-1788-2015}.

%
%
%
%
``{\em The tightest and valid modes apply to all interval types and all operations. The accurate mode is defined only
for inf-sup types}'', that is, to intervals represented by two floating-point numbers of a given format \(\mathbb{F}\) as their endpoints ``{\em (because it involves the nextOut function), and for interval forward and reverse arithmetic
operations \(\mathbf{f}\)}.''
Let \(\mathbf{f}_{\mathrm{exact}}\) denote the corresponding exact mathematical operation, and
\(\mathbf{f}_{\mathrm{tightest}}\) the operation that returns the tightest representable result in the chosen format.
For an input interval \(\mathbf{x}\) representable in the chosen format, 
the definition of the accurate mode ``{\em asserts that \(\mathbf{f}(\mathbf{x})\) is valid, and is at most slightly wider than the result of applying the tightest version to a slightly wider input box:
\[ \mathbf{f}(\mathbf{x}) \subseteq \mathrm{nextOut}(\mathbf{f}_{\mathrm{tightest}}(\mathrm{nextOut}(\mathbf{x}))). \] 
Here the \(\mathrm{nextOut}\) function is defined in terms of the functions \(\mathrm{nextUp}\) and \(\mathrm{nextDown}\) over the set \(\mathbb{F}\) of floating-point numbers as follows.
For \(x \in \mathbb{F}\), define \(\mathrm{nextUp}(x)\) to be \(+\infty\) if \(x = +\infty\), and the least member of \(\mathbb{F}\) greater than \(x\) otherwise; since \(+\infty \in \mathbb{F}\), this is well-defined.
Define \(\mathrm{nextDown}(x)\) to be \(-\mathrm{nextUp}(-x)\).
Then, if \(\mathbf{x}\) is a representable interval, define
\[ \mathrm{nextOut}(\mathbf{x}) = \left\{ \begin{array}{l}
[\mathrm{nextDown} (\underline{x}), \mathrm{nextUp}(\bar{x}) ] \\
\hspace*{12mm} \mbox{ if } \mathbf{x}=[\underline{x},\bar{x}] \neq \emptyset, \\
\emptyset 
\hspace*{10mm} \mbox{ if } \mathbf{x} = \emptyset.
\end{array} \right. \]
For an IEEE 754 format, \(\mathrm{nextUp}\) and \(\mathrm{nextDown}\) are equivalent to the corresponding functions in 5.3.1 of IEEE Std 754-2008.
}


Furthermore, Section 12.10.3, entitled {\em Documentation requirements}, mandates that {\em
An implementation shall document the tightness of each of its interval operations \([\ldots]\)}. 

A last point worth mentioning is that an implementation can provide different accuracy modes.
Here is an excerpt from Section 12.12:
{\em An implementation, or a part thereof, that is 754-conforming shall provide mixed-type operations [\ldots]
for the following operations [\ldots]}
\[\mathrm{add, sub, mul, div, recip, sqrt, sqr, fma.}\]
{\em
An implementation may provide more than one version of some operations for a given type. For instance, it
may provide an ``accurate" version in addition to a required ``tightest" one, to offer a trade-off of accuracy
versus speed or code size. How such a facility is provided is language- or implementation-defined.
}

\section{Devising tests for the accurate mode}
\label{sec:testing}
\subsection{Unit test: definition}

Our focus in this work (in progress) is to test interval arithmetic libraries.
Our approach is based on unit testing. Unit testing consists in testing the implementation of one function or operation (abbreviated as "function" in what follows), independently -- as much as possible -- of the rest of the library. It means completing the following steps:
\begin{enumerate}
\item elaboration of a set of pairs of input(s) and expected output(s), called ``{\em testing pairs}'';
\item evaluation of the function on each input(s);
\item comparison of the result returned by the function with the expected output: usually the test for this pair of input(s) and output(s) is passed if the result equals the expected output, and failed otherwise. More on this later.
\end{enumerate}

Unit testing is commonly used to check whether libraries are either implemented properly, or running correctly on the user's platform.
It is a simple and easy tool to verify the quality of a software development.
Of course it needs to be complemented with other approaches and tools, such as integration testing, that tests whether a whole program computes the expected results, and thus if the functions of the library perform properly together.
Unit testing is in no case a proof that the library is correct; one should resort to formal proof to get such a guarantee, such as those provided by Gappa~\cite[Section 13.3]{HandbookFP} and Flocq~\cite{Flocq-Gappa}.
However, it is simple to implement and very useful, especially at early stages of the software development, to detect bugs rapidly and correct them.
If the set of testing pairs also contains exceptional cases, cornercases, difficult cases, then the number of passed/failed tests is a relevant first approach to asses the quality of the tested library.
This is the only metrics we use: other metrics presented in the literature are not relevant for our work, as our tests are oblivious to the details of the implementation and are directed only towards the mathematical results.

In line with~\cite{PPAM2022-RBFZ}, our focus is to test the functions that are either required or recommended by the IEEE 1788-2015 standard.
For each function \(f\), a testing pair is a pair \( (\mathbf{x}, \mathbf{y})\) where \( \mathbf{x}\) is an input argument and \( \mathbf{y}\) is the tightest representable interval containing \(f(\mathbf{x})\) in the given output format.
A test is passed if \(\mathbf{z}\), the evaluation of \(f(\mathbf{x})\) by the library, verifies \( \mathbf{z} \supseteq \mathbf{y}\): this corresponds to the ``valid mode'' and is the minimal requirement.
Indeed, if \(\mathbf{y}\) were not the tightest representable interval, that is, the minimal interval for inclusion, then a library could compute \(\mathbf{z}\) such that \( \mathbf{z} \subsetneq \mathbf{y}\) and still
\(\mathbf{z} \supseteq f(\mathbf{x})\), that is, \(\mathbf{z}\) satisfies the fundamental principle of interval arithmetic.

\subsection{Testing the ``tightest mode'' and the ``accurate mode''}
\subsubsection*{Testing the ``tightest mode''}
As \(\mathbf{y}\) is the tightest representable result, if the library claims that it provides the tightest mode, then the test is passed if \(\mathbf{z} = \mathbf{y}\).
The question of the determination of \(\mathbf{y}\) can be solved in two ways.
A first, classical, solution is to compute
the endpoints \(\underline{y}\) and \(\bar{y}\) of \(\mathbf{y}=[\underline{y}, \bar{y}]\) using a precision higher than the computing precision used by the library.
Let us illustrate this with a library using {\tt binary64} floating-point numbers, on the (simple: it is monotonic) example of the logarithm function. The given floating-point implementation of \(\log\) does not need to provide correct rounding, however we assume that, for any precision \(q\), it satisfies \(RD_q(\log(x)) \leq \log(x) \leq RU_q (\log (x)) \), for \(RD_q\) rounding downwards, and \(RU_q\) rounding upwards, both in precision \(q\).
In order to compute the infimum \(\underline{y}\) of \(\mathbf{y} = \log \mathbf{x} = \log ([ \underline{x}, \bar{x}])\) where \(\underline{x} >0\), one first computes the approximations \(RD_q(\log(\underline{x}))\) and \( RU_q(\log(\underline{x}))\) of \(\log(\underline{x})\) in high precision \(q\), and finally round them downwards in the target precision, \(53\) for {\tt binary64}.
If \(RD_{53}(RD_q(\log(\underline{x})))\) and \(RD_{53}(RU_q(\log(\underline{x})))\) are equal, then they are the sought value \(\underline{y}\) for the infimum of \(\log(\mathbf{x})\).
Otherwise, a larger value of \(q\) is chosen and this process is repeated.

An easier solution is to assume that the MPFR library~\cite{MPFR}, or the CORE-MATH library~\cite{Core-math2022} once it is fully developed, 
is correct and to use it, e.g. using MPFI~\cite{MPFI}, to compute the reference result \(\mathbf{y}\).
Indeed, depending on the implemented interval arithmetic, a library will need CORE-MATH for binary32 and binary64 floating-point formats, or MPFR for other precisions, in particular arbitrary precision, as in MPFI.
MPFR is a well-know library for arbitrary precision floating-point arithmetic, that offers correct rounding of the result, for any chosen rounding mode and for a large variety of functions that encompasses the set of required and recommended functions of the IEEE 1788-2015 standard.
MPFI is a library for interval arithmetic based on MPFR: it handles intervals, represented by their endpoints which are MPFR numbers, that is, arbitrary precision floating-point numbers.
CORE-MATH is a library under development, its goal is to offer an implementation of the mathematical functions, mainly for {\tt binary32} and {\tt binary64} floating-point numbers, that both offers correct rounding and competitive performance.
The use of a correctly rounded floating-point arithmetic, provided by a library such as CORE-MATH or MPFR, greatly eases the development of an interval arithmetic library. As developers of MPFI, we are well aware that for non-monotonic functions, this development still requires much care in order to determine the tightest result.

\subsubsection*{Testing the ``accurate mode''}
If one relies on MPFR, then one can use it also to determine \(\mathbf{y}'\) the largest (for inclusion) interval that satisfies the ``accurate mode'' condition. 
As MPFR offers the {\tt nextUp} and {\tt nextDown} functions (called {\tt mpfr\_nexabove} and {\tt mpfr\_nextbelow} in MPFR), it is easy to devise the {\tt nextOut} function defined in Section~ \ref{subsec:accurate-mode-def},
to add it to MPFI, and to evaluate \(\mathbf{y}'\) as \(\mathtt{nextOut}(f(\mathtt{nextOut}(\mathbf{x})))\).
A test for the accurate mode is passed if \( \mathbf{y} \subseteq \mathbf{z} \subseteq \mathbf{y}'\).

\subsubsection*{More details about the testing pairs}
We have proposed and provided in a public repository testing pairs for four functions: cubic root, exponential, sine, hyperbolic arc-tangent.
Around a hundred such pairs are currently proposed for each function.

For non-monotonic functions, intervals with endpoints on both sides, and close to an extremum are proposed. 

MPFI is chosen as the reference library to establish the testing pairs. For the "tightest mode", these pairs are trusted to be correct not only because MPFI is trustworthy, but also because they have been used (and checked beforehand) by several interval arithmetic libraries. However, the testing pairs for the accurate mode have been used less broadly.

\subsection{Testing some interval arithmetic libraries}
Few libraries of interval arithmetic are compliant with the IEEE 1788-2015 standard. 
In~\cite{hidden-ref}, we have started to test some of these libraries. We detail below the points relevant for the accurate mode, for the following libraries:
\begin{itemize}
\item the {\tt libieeep1788}~\cite{libieeep1788}, a C++ library, has been developed along the elaboration of the IEEE 1788-2015 standard, as a proof of concept. Operations are implemented using MPFR~\cite{MPFR} for operations on arbitrary-precision floating-point endpoints: thus the library inherits from MPFR the property of returning the tightest results. The accurate mode is not provided, it was probably not specified in the standard when this library was developed. We did not test {\tt libieeep1788} because it is no more maintained and very difficult to install with receent versions of C++.
\item probably one of the oldest compliant library is {\tt JInterval}~\cite{JInterval}, in Java. This library uses either an exact representation of the endpoints (using rational arithmetic), or a representation using the IEEE 754 floating-point format. It claims to offer the tightest mode. Even if we have been able neither to check the code nor to test the library, we are rather confident that this is the case, given the chosen approach, based on exact or arbitrary-precision representation of the endpoints and on the use of the corresponding arithmetic on endpoints.
\item {\tt Octave/Interval}~\cite{OctaveInterval} is probably the easiest-to-use library, as it is based on Octave to benefit from its simplified syntax and its expressiveness. (This comes at the expense of the performance, as the code is interpreted.) Octave offers only the tightest mode, it achieves this accuracy by using the {\tt CRlibm} library~\cite{CRlibm} for mathematical functions.
\item {\tt JuliaIntervals}~\cite{JuliaIntervals} is a recent library, written in the Julia programming language that offers both expressiveness and performance.
For testing purposes, we consider only the underlying floating-point formats for the endpoints that are either the {\tt binary32} or the {\tt binary64} formats of IEEE 754.
{\tt JuliaIntervals} claims to provide the tightest mode, and indeed it successfully passed the tests for most functions.
However, for \( \mathrm{atan}\) and \(\mathrm{cbrt}\), these tests have uncovered that the implementation uses only {\tt binary64} and thus that a double rounding occurs when it is then rounded to {\tt binary32}: thus, for these functions, {\tt JuliaIntervals} currently provides the accurate mode. This observation may not hold for very long, as the authors of {\tt JuliaIntervals} are now aware that the implementation was too naive and they are working on offering tight results for these functions as well.
\end{itemize}

\subsection{Testing the tester}
Our concern regards the verification of the tools used to elaborate the testing pairs.
In this work, the focus is to check, even minimally, whether the {\tt nextOut} function is correct.
To begin with, this check, like our checks for other procedures, relies on unit testing.
We are looking for testing pairs for {\tt nextOut}.
We checked the testing pairs for {\tt nextUp} and {\tt nextDown} (or rather {\tt mpfr\_nextabove} and {\tt mpfr\_nextbelow}) in MPFR.
What we found are, on the one hand, tests for special and exceptional values such as \(0\), \(\infty\) or {\tt NaN} and some powers of \(2\), and on the other hand a random choice of arguments.
For these randomly chosen arguments, the test is passed if, for instance for a value \(x > 0\), both 
{\tt nextUp}\((x)\) and \(x\) plus a small positive quantity round to the same value: this approach is backed up by a mathematical proof that it is correct.

Thus our testing pairs for {\tt nextOut} will include:
\begin{itemize}
\item specific and exceptional values for the endpoints;
\item random values for the endpoints and a test using the same formulas as those used by MPFR;
\item specific cases for intervals: \(\emptyset\), \(\mathbb{R}\), invalid intervals such as \([a,b]\) with \(b < a\);
\item cases to test symmetries: is {\tt nextOut}\(([-b,-a])\) equal to \(-\){\tt nextOut}\([a,b]\) for \(a \in \mathbb{R}\), \(b \in \mathbb{R}\) and \(a \leq b\). Symmetries will also be tested for special and exceptional values.
\end{itemize}
Are these cases sufficient to exhibit one sample for each of the various possibilities?
If not, what is missing? We look forward to get feedback and recommendations from our readers.

\section{Discussion and conclusion}
\label{sec:ccl}
The accurate mode of the IEEE 1788-2015 standard is well defined and the testing methodology well understood and ready to be applied.
However, this accurate mode has not been implemented so far, and the testing tool did not have a chance to be employed.
The development of this tool has nonetheless been a good opportunity to dig deeper into the notion of {\em accurate mode}.

We now have some concerns about the pertinence and the viability of this notion of accurate mode:
\begin{itemize}
\item no IEEE 1788-2015 compliant library implements it, the preference goes to the tightest mode;
\item with the MPFR library for arbitrary precision floating-point numbers and the CORE-MATH library for the {\tt binary32} and {\tt binary64} floating-point numbers, there will be no gain in terms of performance between the implementations of the tightest mode and of the accurate mode;
\item there is a surprising phenomenon that results from the definition of the accurate mode: in {\tt binary64}, the largest accurate evaluation of 
\(\sin ([0,10])\) is \([-1.0000000000000003e0,1.0000000000000003e0]\)
and of \(\exp ([-10^9,0])\) is \([-4.9406564584124655e-324,1.0000000000000005e0]\),
when both are printed in decimal and rounded away.
A result strictly enclosing \([-1,1]\) for the sine function is unexpected and even disturbing. So is a strictly negative result for the exponential function.
A mathematically more satisfying definition would prohibit results so clearly out of the well-known range for these functions.
This could take the form of a long list of exceptions to the general definition: this will be suggested to the revision committee for IEEE 1788-2015.
\end{itemize}

A provocative conclusion could be: shall the IEEE 1788-2015 standard abandon the notion of accurate mode?
As the revision of the IEEE 1788-2015 must be performed before the end of 2025, these remarks
are intended as food for thought and will be brought to the attention of the committee.

\end{document}